\input amstex 
\documentstyle{amsppt}
\input bull-ppt
\keyedby{bull382/pxp}

\topmatter
\cvol{28}
\cvolyear{1993}
\cmonth{April}
\cyear{1993}
\cvolno{2}
\cpgs{329-333}
\title On the distribution of sums of residues \endtitle
\author Jerrold R. Griggs\endauthor
\shortauthor{J. R. Griggs}
\shorttitle{Sums of residues}
\address Department of Mathematics, University of South 
Carolina,  Columbia,
South Carolina 29208\endaddress
\ml griggs\@math.scarolina.edu\endml
\date September 2, 1992\enddate
\thanks Research supported in part by NSA/MSP Grant 
MDA90-H-4028 and by a
Visiting Professorship at Simon Fraser University\endthanks
\abstract We generalize and solve the $\roman{mod}\,q$ 
analogue of a problem of
Littlewood and Offord, raised by Vaughan and Wooley, 
concerning the
distribution of the $2^n$ sums of the form 
$\sum_{i=1}^n\varepsilon_ia_i$,
where each $\varepsilon_i$ is $0$ or $1$. For all $q$, 
$n$, $k$ we determine
the maximum, over all reduced residues $a_i$ and all sets 
$P$ consisting of $k$
arbitrary residues, of the number of these sums that belong 
to $P$.\endabstract
\endtopmatter

\document
\heading 1. Introduction\endheading
Vaughan and Wooley \cite{15} raised the problem of 
determining the maximum
number of the $2^n$ sums of the form 
$\sum_{i=1}^n\varepsilon_ia_i$, where each
$\varepsilon_i$ is $0$ or $1$, that are congruent to 
$0\mod q$. The maximum is
over all residues $a_1,\dotsc,a_n$ that are {\it 
reduced\/}, which means that
$(a_i,q)=1$ for all $i$. Results about this problem have 
been applied to study
the solutions of simultaneous additive equations.

By using analytical tools, including exponential sums and 
classical
inequalities, and by treating many cases, Vaughan and 
Wooley show that the
maximum is $\(\smallmatrix n\\ \lfloor 
n/2\rfloor\endsmallmatrix\)$ provided
that $q>\lceil n/2\rceil$. This bound is sharp, since it 
is attained by 
letting $a_i$ be $1$ for $i\le \lceil n/2\rceil$ and 
$-1$ for $i>\lceil n/2\rceil$. (To see this, observe that 
for
this choice of $a_i$'s, we have 
$\sum_{i=1}^n\varepsilon_ia_i\equiv 0$
precisely when an equal number of $\varepsilon_i$'s are 
$1$ for 
$i\le\lceil n/2\rceil$ and $-1$ for $i>\lceil n/2\rceil$. 
This happens if and only if the number of indices $i$ with
$i\le \lceil n/2\rceil$ and $\varepsilon_i=1$ plus the 
number with $i>\lceil
n/2\rceil$ and $\varepsilon_i=0$ is $\lfloor n/2\rfloor$, 
so that the choices
correspond to the subsets of $\{1,\dotsc,n\}$ of size 
$\lfloor n/2\rfloor$.)

When $\lceil n/2\rceil\ge q$, wraparound effects $\mod q$ 
come into play. For
example, with the $a_i$'s chosen as above, the sum 
$\sum_{i=1}^qa_i$ is also
congruent to $0$, so the answer exceeds $\(\smallmatrix 
n\\ \lfloor
n/2\rfloor\endsmallmatrix\)$. 

We solve the problem for arbitrary $n$ and $q$, using an 
inductive argument
that is inspired by the study of the extremal properties 
of the Boolean lattice
$B_n$ on the collection $2^{[n]}$ of all subsets of the 
$n$-set
$[n]=\{1,\dotsc,n\}$, ordered by inclusion. Let us adopt 
the notation
$$
\binom ns_q\coloneq|\{A\subseteq[n]\:|A|\equiv 
s\}|=\sum_{j\equiv s}\binom js,
$$
for the $\roman{mod}\, 
q$ binomial coefficients in $n$. We shall see that for 
general
$n$ and $q$, the maximum number of sums congruent to $0$ 
is the middle $\mod q$
binomial coefficient $\(\smallmatrix n\\ \lfloor 
n/2\rfloor\endsmallmatrix\)_q$.
The maximum is attained as before by dividing the $a_i$'s 
as evenly as possible
between $1$ and $-1$. In general, the maximum number of 
sums congruent to any
single residue is $\(\smallmatrix n\\ \lfloor 
n/2\rfloor\endsmallmatrix\)_q$.
Throughout the paper we maintain the condition that the 
residues $a_i$ be
reduced. Without this restriction, one would select 
$a_i$'s with common
factors,
in order to increase the number of sums congruent to $0$.

This problem we are considering is the analogue for 
residues of a famous
problem about the clustering of partial sums of a 
collection of complex
numbers. In connection with their study of roots of random 
polynomials,
Littlewood and Offord \cite{13} were led to consider the 
following question.
For $a_1,\dotsc,a_n\in\bold C$ with $\|a_i\|\ge 1$ for all 
$i$ and for an open
ball $S\subset\bold C$ of unit diameter, how many of the 
$2^n$ sums of the form
$\sum_{i=1}^n\varepsilon_ia_i$, where each $\varepsilon_i$ 
is $0$ or $1$, can
belong to $S$? They sought the maximum over all choices of 
$a_i$'s and $S$.
In particular, if one selects $a_i$ to be $1$ for all $i$ 
and centers $S$ at 
$\lfloor n/2\rfloor$, one can pack $\(\smallmatrix n\\ 
\lfloor
n/2\rfloor\endsmallmatrix\)$ sums into $S$, and this was 
believed to be
optimal. While Erd\"os \cite3 soon proved this for the 
real $a_i$ case, it was
twenty years before the original complex case was solved, 
by Katona \cite8 and
Kleitman \cite{10} independently, from an appropriate 
extension of the theorem
of Sperner \cite{14} about the maximum size of an 
antichain in the Boolean
lattice.

Although the usual Sperner method does not extend to 
higher dimensions,
Kleitman \cite{11, also in 5} found a remarkable proof 
that the answer is still 
$\(\smallmatrix n\\ \lfloor n/2\rfloor\endsmallmatrix\)$ 
in any dimension $m$
(or indeed, in Hilbert space): For any vectors 
$a_1,\dotsc, a_n\in\bold R^m$ of
length at least one, there is a partition of $2^{[n]}$ 
into just 
$\(\smallmatrix n\\ \lfloor n/2\rfloor\endsmallmatrix\)$ 
blocks, such that for
any sets $I$, $J$ in the same block, the sums $\sum_{i\in 
I}a_i$ and
$\sum_{j\in J}a_j$ are far apart (distance at least one). 
Hence any open ball
$S$ of unit diameter contains at most $\(\smallmatrix n\\ 
\lfloor
n/2\rfloor\endsmallmatrix \)$ sums. The idea behind this 
construction is that
for every $n$, the sizes of the blocks partitioning 
$2^{[n]}$ exactly match the
sizes of the chains in the famous inductive symmetric 
chain decomposition of
$B_n$ discovered by de Bruijn et al. \cite2.

Erd\"os \cite3 considered the more general problem of 
maximizing the number of
sums of vectors inside an open ball  in $\bold R^m$ of 
diameter $d\ge 1$. He
solved this problem for the real case $(m=1)$ using 
Sperner theory, and he
found that the value attained when all $a_i=1$ is optimal 
for all $d$. This
value is the sum of the $\lceil d\rceil$ 
middle binomial coefficients, $\sum_{(n-\lceil
d\rceil)/2\le j<(n+\lceil d\rceil)/2}\binom nj$. However, 
the problem is more
complicated when $m\ge 2$ and completely solved only in 
some special cases.
Using a variety of tools from extremal set theory, 
probability, and geometry,
many authors have attacked this more general question, 
including Kleitman
\cite{12}, Griggs \cite6, and Frankl and F\"uredi \cite4. 
Also see the survey
by Anderson \cite1.

In marked contrast to previous results of the 
Littlewood-Offord type, the
setting for the work of Vaughan and Wooley is the additive 
group $\bold Z_q$ of
integers $\mod q$. Nonetheless, as with the unit diameter 
problem above, we
shall see that their theorem can be obtained by an 
inductive partition
construction inspired by a particular chain partition of 
the Boolean lattice.
The method yields the solution to the extension of their 
problem to general $n$
and $q$.

More generally, we determine the maximum number of the 
$2^n$ sums
$\sum_{i=1}^n
\varepsilon_ia_i$ congruent $\mod q$ to any of $k$ 
arbitrary residues 
$\rho_j$, for $1\le j\le k$, over all choices of the 
residues $\rho_j$ and the
reduced residues $a_i$. The answer is the sum of the $k$ 
middle $\mod q$
binomial coefficients in $n$. This bound is attained by 
selecting all $a_i$ to
be $1$ and selecting the $k$ middle values for the 
residues $\rho_j$. Switching
some $a_i$ in this solution to $-1$ has the effect of 
shifting the collection
of all $2^n$ sums down by $1$. Thus the bound is also 
attained by selecting
$a_i$ to be $1$ for $i\le\lceil n/2\rceil$ and $-1$ for 
$i>\lceil n/2\rceil$
and by choosing the $k$ initial values in the sequence 
$0,1,-1,2,-2,\dotsc$ for
the residues $\rho_j$.

\heading 2. The main result\endheading
We fix the integer $q>0$ and work in $\bold Z_q$.

\thm{Theorem 1} Let $a_1,\dotsc,a_n$ be reduced residues 
in $\bold Z_q$. Let
$P\subseteq\bold Z_q$, where $|P|=k$. Then the number of 
the $2^n$ sums
$\sum_{i=1}^n\varepsilon_ia_i$ in $P$, where each 
$\varepsilon_i$ is $0$ or
$1$, is at most the sum of the $k$ middle $\mod q$ 
binomial coefficients
$\sum_{(n-k)/2\le j<(n+k)/2}\binom nj_q$, and this bound 
is best possible.
\ethm

\demo{Proof} For $S\subseteq \bold Z_q$ and $a\in\bold 
Z_q$, let
$S+a\coloneq\{s+a\:s\in S\}\subseteq\bold Z_q$. For $\scr 
A\subseteq 2^{[n]}$,
define the {\it sum set\/}
$$
S(\scr A)=\left\{\sum_{i\in I}a_i\pmod{q}\:I\i\scr 
A\right\}\.
$$
We say that $\scr A\subseteq 2^{[n]}$ is a {\it 
structure\/} for
$a_1,\dotsc,a_n$ provided that the sums in $S(\scr A)$ are 
distinct.

We shall partition $2^{[n]}$ into $\(\smallmatrix n\\ 
\lfloor
n/2\rfloor\endsmallmatrix\)_q$ structures in such a way 
that the bound in the
theorem will follow for all $k$. The construction is 
carried out by induction
on $n$ for a given sequence of reduced residues 
$a_1,a_2,\dotsc$. It starts at
$n=0$ with the single structure $\{\emptyset\}$. For the 
induction step,
suppose we are given a partition of $2^{[n-1]}$ into 
structures $\scr A_j$ for
$a_1,\dotsc,a_{n-1}$. Then the structures $\scr A_j$ and 
$\scr
A_j'\coloneq\{I\cup\{n\}\:I\in \scr A_j\}$ for 
$a_1,\dotsc,a_n$ partition
$2^{[n]}$, but they are not quite the ones we want. Notice 
that 
$S(\scr A_j')=S(\scr A_j)+a_n$. We require an easy fact.

\thm{Lemma} Let $\emptyset\ne S\subseteq\bold Z_q$ and 
$a\in\bold Z_q$ with
$(a,q)=1$. Then $S+a=S$ if and only if $S=\bold Z_q$.
\ethm

If $S(\scr A_j)$ is $\bold Z_q$, then so is $S(\scr 
A_j')$, and we leave both
structures alone. However, if $S(\scr A_j)\ne\bold Z_q$, 
then by the lemma
there exists at least one element $t\in S(\scr 
A_j')\backslash 
S(\scr A_j)$, say 
$t=\sum_{i\in I}a_i$ where $I\in \scr A_j'$, so that we 
may replace $\scr A_j$
and $\scr A_j'$ by the structures $\scr B_j=\scr 
A_j\cup\{I\}$ and $\scr B_j'=
\scr A_j'\backslash\{I\}$. We have $|\scr B_j|=|\scr A_j|+
1$ and $|\scr B_j'|=
|\scr A_j|-1$. In the case where $|\scr A_j|=1$, we 
discard $\scr B_j'$.

Now denote the structures in this partition of $2^{[n]}$ 
by $\scr A_j$ for
$j=1,2,\dotsc$. Since sets in a structure have distinct 
sums, it follows that
$$
\biggl|\biggl\{I\subseteq[n]\:\sum_{i\in I}a_i\in 
P\biggr\}\biggr|\le\sum_j
\min(k,|\scr A_j|)\.\tag 1
$$
It suffices to show that the sum on the right-hand side of 
inequality (1) is at
most the sum of the $k$ middle $\roman{mod}\,q$ binomial 
coefficients in $n$.

Since the collection of structure sizes $|\scr A_j|$ 
depends in no way on the
actual values of the $a_i$'s, it is enough to consider the 
case where all
$a_i=1$. One can verify by induction on $n$ that the sum 
set $S(\scr A_j)$ for
each structure consists of all $q$ residues or else 
consists of values
congruent to an interval $x,
x+1,\dotsc,y\in\bold Z$ centered about $n/2$, which means 
$x+y=n$. The number of structures in the partition is 
$\(\smallmatrix
n\\ \lfloor n/2\rfloor\endsmallmatrix\)_q$, because every 
structure contains a
set with sum (i.e., cardinality) $\equiv\lfloor 
n/2\rfloor$. For general $k$,
we see that the sum on the right-hand side of (1) is the 
number of subsets of
cardinality congruent to any of the $k$ middle values 
around $n/2$.\qed
\enddemo

When $q>\lceil n/2\rceil$, we have that $\(\smallmatrix 
n\\ \lfloor
n/2\rfloor\endsmallmatrix\)_q=\(\smallmatrix n\\ \lfloor
n/2\rfloor\endsmallmatrix\)$, which implies the original 
result of Vaughan and
Wooley \cite{15}:

\thm{Corollary 1} Let $a_1,\dotsc,a_n$ be reduced residues 
in $\bold Z_q$,
where $q>\lceil n/2\rceil$. Then the number of the $2^n$ 
sums $\sum_{i=1}^n
\varepsilon_ia_i$ congruent to $0$, where each 
$\varepsilon_i$ is $0$ or $1$,
is at most $\(\smallmatrix n\\ \lfloor 
n/2\rfloor\endsmallmatrix\)$, and this
bound is best possible.
\ethm

\heading 3. Related remarks\endheading
The inspiration for the proof of the theorem is the 
inductive partition of the
Boolean lattice $B_n$ into saturated chains of size at 
most $q$, that is, into
collections of at most $q$ totally ordered subsets of 
consecutive sizes. Katona
\cite9 used this construction to determine the maximum 
number of subsets of 
$\{1,\dotsc,n\}$ containing no sets $A\subset B$ with 
$0<|B\backslash A|<q$.
The author \cite7 later independently devised the same 
construction to obtain a
maximum-sized collection of disjoint saturated chains of 
size $q$ in $B_n$. The
collection of structure sizes $|\scr A_j|$ in our 
construction exactly
corresponds to the collection of chain sizes in Katona's 
partition.

By applying the theorem with $k=q-1$, it is also possible 
to determine the
minimum number of sums in any residue class.

\thm{Corollary 2} Let $a_1,\dotsc,a_n$ be reduced residues 
in $\bold Z_q$,
where $n\ge q-1$. Let $\rho\in\bold Z_q$. Then the number 
of the $2^n$ sums
$\sum_{i=1}^n\varepsilon_ia_i$ congruent to $\rho$, where 
each $\varepsilon_i$
is $0$ or $1$, is at least $\(\smallmatrix n\\
\lceil(n-q)/2\rceil\endsmallmatrix\)_q$, and this bound is 
best possible.
\ethm

The bound in Corollary 2 is attained by taking all $a_i=1$ 
and
$\rho\equiv\lceil(n-q)/2\rceil$. For $n<q-1$, no sums are 
congruent to $-1$
when all $a_i$ equal $1$. The asymptotic growth of the 
$\roman{mod}\,q$ binomial
coefficients, studied in connection with saturated chain 
partitions \cite7,
implies that the lower bound in Corollary 2 approaches 
$2^n/q$ as $n\to\infty$
with fixed $q$. (This remains true even if $q$ grows with 
$n$, provided that
$q=o(n^{1/2})$.) Hence, for any sequence 
$\{a_1,a_2,\dotsc\}$ of reduced 
residues $\roman{mod}\,q$, the distribution of the 
$\roman{mod}\,q$ sums of the first $n$
residues is asymptotically uniform as $n\to\infty$.

The Littlewood-Offord problem has an equivalent 
formulation that considers the
concentration of sums of the form 
$\sum_{i=1}^n\delta_ia_i$ with each
$\delta_i=1$ or $-1$, where as before $\|a_i\|\ge 1$ for 
all $i$. The analogous
problem in $\bold Z_q$ can be solved by a reduction to the 
original problem of
Theorem 1.

\thm{Corollary 3} Let $a_1,\dotsc,a_n$ be reduced residues 
in $\bold Z_q$.
Let $P\subseteq \bold Z_q$, where $|P|=k$. Then the number 
of the $2^n$ sums
$\sum_{i=1}^n\delta_ia_i$ in $P$, where each $\delta_i$ is 
$1$ or $-1$, is at
most the sum of the $k$ middle $\roman{mod}\,r$ 
binomial coefficients in $n$, where $r$
is $q$ when $q$ is odd and $q/2$ when $q$ is even, and 
this bound is best
possible.
\ethm

\heading Acknowledgment\endheading
The author is grateful to Oren Patashnik and Ted Sweetser 
for many suggestions
that greatly improved the presentation of this paper.

\enddocument